\documentclass[english,12 pt]{article}
\usepackage[T1]{fontenc}
\usepackage[latin1]{inputenc}
\usepackage{amsmath}

\makeatletter

\usepackage{amssymb,amsthm}
\usepackage{enumerate}

\makeatletter
    \newcommand{\figcaption}{\def\@captype{figure}\caption}
    \newcommand{\tabcaption}{\def\@captype{table}\caption}
\makeatother

\usepackage{babel}
\makeatother
\begin{document}

\title{Expected Utility Optimization\\ Calculus of Variations Approach}
\author{Khoa Tran}
\date{\today}
\maketitle
\abstract{In this paper, I'll derive the Hamilton-Jacobi (HJ) equation for Merton's problem in Utility Optimization Theory using a Calculus of Variations (CoV) Approach. For stochastic control problems, Dynamic Programming (DP) has been used as a standard method. To the best of my knowledge, no one has used CoV for this problem. In addition, while the DP approach cannot guarantee that the optimum satisfies the HJ equation, the CoV approach does.}
\section{Introduction}
Consider Merton's problem in Utility Optimization Theory. That is: given a stock with dynamics
\[
dS_t=\mu S_tdt+\sigma S_tdW_t
\]
(where $W_t$ denotes the Brownian motion), find the amount $\Pi_t$ to invest in this stock that maximizes the expected utility at a terminal time $T$. Mathematically, if $X_t$ represents the wealth, the problem to solve is
\[
\max_{\Pi_t}\mathbb{E}[U(X_T|X_0],
\]
where $U$ is the utility function ($U$ is given and is concave).

Suppose that the control function $\Pi_t=\alpha(t,X_t)$ depends only on the current state but not the past\footnote{Is it a reasonable assumption if the stock price follows a Brownian motion, i.e. it has Markov property?}. Then the wealth has the following dynamics\footnote{For simplicity, let's assume that $\mu$ and $\sigma$ are independent of the stock price.}:
\begin{equation}\begin{split}
dX_t &= rX_tdt + (\mu-r)\Pi_tdt +\sigma\Pi_tdW_t\\
&= a(t,X_t)dt+b(t,X_t)dW_t,
\end{split}\end{equation}
where
\begin{align}
a(t,X_t)&=rX_t+(\mu-r)\alpha(t,X_t)\\
b(t,X_t)&=\sigma \alpha(t,X_t).
\end{align}

The Fokker-Planck (FP) equation, which describes the evolution of probability density $p(t,x)$ of the random process $X_t$, has the form
\begin{equation}
\frac{\partial p}{\partial t}=-\frac{\partial (ap)}{\partial x}+\frac{1}{2}\frac{\partial^2(b^2p)}{\partial x^2},
\end{equation}
with the initial condition $p(0,x)=p_0(x)=\delta(x-X_0)$.

Our goal is to maximize $\mathbb{E}[U(X_T)|X_0]=\int_\mathbb{R}U(x)p(T,x)dx$ over all admissible control functions $\alpha(t,X_t)$. Since the utility expectation is a \underline{concave}\footnote{Underlined statement can be proved but is not shown in this paper.} functional of the control, the optimum can be determined by the CoV method. This is actually a PDE-constrained\footnote{the density function is constrained by the FP equation} optimization.
\section{Calculus of Variations}
Let's introduce the Lagrange multipliers $\lambda(t,x)$ and $\lambda_0(x)$ and consider the problem of optimizing
\[\begin{split}
L(p,\alpha,\lambda,\lambda_0)=&\int_\mathbb{R}U(x)p(T,x)dx + \int_\mathbb{R}\lambda_0(x)[p(0,x)-p_0(x)]dx\\
&+\int_0^T\int_\mathbb{R}\lambda\left[\frac{\partial}{\partial t}p+\frac{\partial}{\partial x}(ap)-\frac{1}{2}\frac{\partial^2}{\partial x^2}(b^2p)\right]dxdt.
\end{split}\]
\begin{enumerate}
\item Taking the first variation of $L$ with respect to $p$, we obtain
\[\begin{split}
\delta_pL=&\int_\mathbb{R}U(x)\delta p(T,x)dx + \int_\mathbb{R}\lambda_0(x)\delta p(0,x))dx\\
&+\int_0^T\int_\mathbb{R}\lambda\left[\frac{\partial}{\partial t}\delta p+\frac{\partial}{\partial x}(a\delta p)-\frac{1}{2}\frac{\partial^2}{\partial x^2}(b^2\delta p)\right]dxdt.
\end{split}\]
Note, by integration by parts, that
\[\begin{split}
\int_0^T\int_\mathbb{R}\lambda\frac{\partial}{\partial t}\delta p dxdt =& (\int_\mathbb{R}\lambda\delta pdx)\big |_0^T-\int_0^T\int_\mathbb{R}\lambda_t\delta pdxdt\\
=&\int_\mathbb{R}\lambda(T,x)\delta p(T,x)dx-\int_\mathbb{R}\lambda(0,x)\delta p(0,x)dx\\
&-\int_0^T\int_\mathbb{R}\lambda_t\delta pdxdt,
\end{split}\]
\[\begin{split}
\int_0^T\lambda\frac{\partial}{\partial x}(a\delta p)dxdt =& \int_0^T\lambda a\delta p\big |_{-\infty}^{+\infty}dt -\int_0^T\int_\mathbb{R}\lambda_xa\delta pdxdt\\
=&\int_0^T\lambda(t,+\infty) a(t,+\infty)\delta p(t,+\infty)dt\\
&- \int_0^T\lambda(t,-\infty) a(t,-\infty)\delta p(t,-\infty)dt\\
&-\int_0^T\int_\mathbb{R}\lambda_xa\delta pdxdt,
\end{split}\]
and
\[\begin{split}
\int_0^T\int_\mathbb{R}\frac{\lambda}{2}\frac{\partial^2}{\partial x^2}(b^2p)dxdt =& \int_0^T\frac{\lambda}{2}\frac{\partial}{\partial x}(b^2\delta p)\big |_{-\infty}^{+\infty}dt-\int_0^T\frac{\lambda_x}{2}b^2\delta p\big |_{-\infty}^{+\infty}dt\\
&+\int_0^T\int_\mathbb{R}\frac{\lambda_{xx}}{2}b^2\delta pdxdt.
\end{split}\]
If $\delta_pL=0$ for all admissible $\delta p$ then
\begin{align*}
\lambda(T,x)=&-U(x)\\
\lambda_0(x)=&\lambda(0,x)\\
\lambda_t+a\lambda_x+\frac{b^2}{2}\lambda_{xx}=&0,
\end{align*}
and some boundary (at infinity) conditions must be satisfied.
\item Taking the first variation of $L$ with respect to $\alpha$ (which occurs only in $a$ and $b$), we obtain
\[\begin{split}
\delta_\alpha L=&\int_0^T\int_\mathbb{R}\lambda\frac{\partial}{\partial x}[(\mu-r)p\delta\alpha]dxdt-\int_0^T\int_\mathbb{R}\lambda\frac{\partial^2}{\partial x^2}(\sigma^2p\alpha\delta\alpha)dxdt\\
=&\int_0^T\lambda(\mu-r)p\delta\alpha\big |_{-\infty}^{+\infty}dt-\int_0^T\int_\mathbb{R}\lambda_x(\mu-r)p\delta\alpha dxdt\\
&-\int_0^T\lambda\frac{\partial}{\partial x}(\sigma^2p\alpha\delta\alpha)\big |_{-\infty}^{+\infty}+\int_0^T\int_\mathbb{R}\lambda_x\frac{\partial}{\partial x}(\sigma^2p\alpha\delta\alpha)dxdt\\
=&\int_0^T\lambda(\mu-r)p\delta\alpha\big |_{-\infty}^{+\infty}dt-\int_0^T\int_\mathbb{R}\lambda_x(\mu-r)p\delta\alpha dxdt\\
&-\int_0^T\lambda\frac{\partial}{\partial x}(\sigma^2p\alpha)\delta\alpha\big |_{-\infty}^{+\infty}dt - \int_0^T\lambda\sigma^2p\alpha\frac{\partial}{\partial x}\delta\alpha\big |_{-\infty}^{+\infty}dt\\
&+\int_0^T\lambda_x\sigma^2p\delta\alpha\big|_{-\infty}^{+\infty}dt-\int_0^T\int_\mathbb{R}\lambda_{xx}\sigma^2p\alpha\delta\alpha dxdt.
\end{split}\]
If $\delta_\alpha L=0$ for all admissible $\alpha$ then
\[
\lambda_{xx}\sigma^2p\alpha+\lambda_x(\mu-r)p=0,
\]
or $\sigma^2\lambda_{xx}\alpha+(\mu-r)\lambda_x=0$ almost everywhere \underline{for $p\neq 0$ a.e.}) and some boundary conditions must be satisfied.
\item Setting $\delta_\lambda L=0$ and $\delta_{\lambda_0}L=0$ gives us the Fokker-Planck equation with the initial condition of $p$.
\end{enumerate}
\underline{Since $L$ is concave}, the necessary and sufficient condition for optimality can be described by the following equations (with boundary conditions at infinity not stated here).
\begin{itemize}
\item State equations:
\begin{align}
\frac{\partial p}{\partial t}=&-\frac{\partial (ap)}{\partial x}+\frac{1}{2}\frac{\partial^2(b^2p)}{\partial x^2},\\
p(0,x)=&\delta(x-X_0).
\end{align}
\item Adjoint equations:
\begin{align}
\lambda(T,x)=&-U(x),\\
\lambda_t+a\lambda_x+\frac{b^2}{2}\lambda_{xx}=&0.
\end{align}
\item Control equation:
\begin{equation}
\sigma^2\lambda_{xx}\alpha+(\mu-r)\lambda_x=0.
\end{equation}
\end{itemize}
The control equation implies that
\begin{equation}
\alpha=-\frac{(\mu-r)\lambda_x}{\sigma^2\lambda_{xx}}
\end{equation}
and hence
\begin{align*}
a(t,x)=&rx-\frac{(\mu-r)^2\lambda_x}{\sigma^2\lambda_{xx}},\\
b(t,x)=&-\frac{(\mu-r)\lambda_x}{\sigma\lambda_{xx}}.
\end{align*}
Plugging these into the adjoint equation, we obtain a Hamilton-Jacobi (HJ) equation of $\lambda$:
\begin{equation}
\lambda_t+rx\lambda_x-\frac{(\mu-r)^2\lambda_x^2}{2\sigma^2\lambda_{xx}}=0.
\end{equation}
The Lagrange multiplier $\lambda$ is essentially the same as the function value $V$ in the Dynamic Programming (DP) approach\footnote{In DP approach, the value function $V(t,x)$ denotes the optimized expected utility at $T$ given initial wealth $x$ at time $t$.}.
\section{Discussions}
\begin{itemize}
\item For variable drift and volatility (i.e. they depend on the stock price), the approach can be extended easily by using the two dimensional version of FP equation (actual derivation will be shown later). For stochastic volatility? I'm not sure but will think more about it.
\item The method can also be extended easily if we want to optimize the expected utility over a period of time instead of at the terminal time. But it's likely that the result can also be derived easily by means of DP.
\item In DP approach, if we can solve the derived HJ equation then we know that it's the optimum. However, it doesn't guarantee that the optimum has to satisfy the HJ equation\footnote{If this sounds mysterious, please refer "Robert C. Merton, 1973. \itshape{An Intertemporal Capital Asset Pricing Model}. Econometrica 41: 867-887".}. The CoV approach shows that $\lambda$ has to satisfy the HJ equation\footnote{It has been shown Crandall and Lions that for HJ-type equations, the viscosity solution, which is the appropriate solution in this context, exists uniquely.}.
\item How do we explain that the Lagrange multiplier $\lambda$ is the value function? This fact seems to have been well explained by economists. (Baxley and Moohouse, `Lagrange multiplier problems in economics', The Am. Math. Monthly.)
\end{itemize}
\section{Numerical methods for computing viscosity solution}
This is what I currently want to study.
%
\section{Derivation of the Fokker-Planck equation}
In this section, we'll derive the one-dimensional Fokker-Planck equation that describes the evolution of density of the stochastic variable $X_t$ with dynamics
\[
dX_t=\mu(t,X_t)dt + \sigma(t,X_t)dW_t.
\]

For any function $f(x)$, consider $\langle f(X_t)\rangle=\frac{d}{dt}\mathbb{E}[f(X_t)]$ at an arbitrary time $t$. On one hand,
\[
\frac{d}{dt}\langle f(X_t\rangle=\frac{d}{dt}\int_\mathbb{R}p(t,x)f(x)dx=\int_\mathbb{R}p_t(t,x)f(x)dx.
\]
On the other hand,
\[\begin{split}
\frac{d}{dt}\langle f(X_t)\rangle=&\frac{\langle df(X_t)\rangle}{dt}\\
=&\langle\mu(t,X_t)f_x+\frac{1}{2}\sigma^2(t,X_t)f_{xx}\rangle \hspace{15pt}\text{(by Ito's formula)}\\
=&\int_\mathbb{R}p(t,x)\left[\mu(t,x)f_x(x)+\frac{1}{2}\sigma^2(t,x)f_{xx}(x)\right]dx\\
=&\int_\mathbb{R}f(x)\left[-\frac{\partial}{\partial x}(\mu p)+\frac{1}{2}\frac{\partial^2}{\partial x^2}(\sigma^2p)\right]dx\hspace{15pt}\text{(integration by parts)}.
\end{split}\]
Since the two formulations hold for any function $f$,
\begin{equation}
p_t=-\frac{\partial}{\partial x}(\mu p)+\frac{1}{2}\frac{\partial^2}{\partial x^2}(\sigma^2p).
\end{equation}
\end{document}